\documentclass[12pt]{amsart} 

\usepackage{amsmath,amsthm,amsfonts,latexsym,euscript,amssymb}

\newtheorem{DF}{Definition}[section]
\newtheorem{LM}[DF]{Lemma}
\newtheorem{PROP}[DF]{Proposition}
\newtheorem{THM}[DF]{Theorem}
\newtheorem{COR}[DF]{Corollary}
\newtheorem{RMK}[DF]{Remark}
\newtheorem{RMKS}[DF]{Remarks}
\newtheorem{PROB}[DF]{Problem}

\newcommand{\bgdf}{\begin{DF}}
\newcommand{\nddf}{\end{DF}}
\newcommand{\bglm}{\begin{LM}}
\newcommand{\ndlm}{\end{LM}}
\newcommand{\bgprop}{\begin{PROP}}
\newcommand{\ndprop}{\end{PROP}}
\newcommand{\bgth}{\begin{THM}}
\newcommand{\ndth}{\end{THM}}
\newcommand{\bgthm}{\begin{THM}}
\newcommand{\ndthm}{\end{THM}}
\newcommand{\bgcor}{\begin{COR}}
\newcommand{\ndcor}{\end{COR}}
\newcommand{\bgrmk}{\begin{RMK}}
\newcommand{\ndrmk}{\end{RMK}}
\newcommand{\bgrmks}{\begin{RMKS}}
\newcommand{\ndrmks}{\end{RMKS}}
\newcommand{\bgprob}{\begin{PROB}}
\newcommand{\ndprob}{\end{PROB}}
\newcommand{\bgeq}{\begin{eqnarray}}
\newcommand{\ndeq}{\end{eqnarray}}
\newcommand{\bgeqq}{\begin{eqnarray*}}
\newcommand{\ndeqq}{\end{eqnarray*}}

\numberwithin{equation}{section}

\newcommand{\dfref}[1]{Definition~\ref{#1}}

\newcommand{\propref}[1]{Proposition~\ref{#1}}

\newcommand{\thmref}[1]{Theorem~\ref{#1}}

\newcommand{\probref}[1]{Problem~\ref{#1}}
\newcommand{\secref}[1]{\S\ref{#1}}
\begin{document}



\title[Classifying Simple Quantum Groups]
{\bf
On the Problem of Classifying
\\ Simple Compact Quantum Groups
}
\author[Shuzhou Wang]
{\bf Shuzhou Wang}
\thanks{}


\address{Department of Mathematics, University of Georgia,
Athens, GA 30602, USA
\newline \indent
Fax: 706-542-5907,
Tel: 706-542-0884
}
\email{szwang@math.uga.edu}



\dedicatory{Dedicated to Professor S. L. Woronowicz on the Occasion of his 70th Birthday}

\begin{abstract}
We review the notion of simple compact quantum groups and
examples, and discuss the problem of construction and
classification of simple compact quantum groups. 
Several new quantum groups constructed by 
Banica, Curran and Speicher since the author's first paper 
on simple quantum groups are shown to be simple using 
results of Raum and Weber.  
\end{abstract}

\maketitle

\section{Introduction}
\label{Introd}

There are two main areas in the operator algebraic approach to quantum groups:
compact quantum groups and locally compact quantum groups.
The former has a satisfactory axiomatic theory due to Woronowicz \cite{Wor87b,Wor98a}.
While the latter has witnessed great progress
due to concerted efforts of several generations of mathematicians
(cf. an incomplete list including
\cite{Kac63a65a,KV74a,EnSch75a,EnSch92a,BS93,KusVaes00a,MaNaWor03a}),
it still does not have an axiomatic framework
that contains the non-compact Drinfeld-Jimbo quantum groups
\cite{Dr87a} as examples except for special cases such as
$SL_q(2, {\mathbb C})$ (cf. e.g. \cite{PW90}),
nor has the existence of Haar weight been established in general.
This is in stark contrast with the fact that
compact real forms of the Drinfeld-Jimbo quantum groups are special
examples in the theory of compact quantum groups
(see \cite{Rosso87a,Rosso90a,SV,Soib1,Lev1}), and Haar measure,  
Peter-Weyl theorem and Tannaka-Krein duality can be established 
from a very simple set of axioms \cite{Wor87b,Wor88a,Wor98a}. 
Within this well established framework of compact quantum groups, 
recent work on compact quantum groups has been primarily on 
the construction, classification, structure and other (operator) algebraic 
properties of specific classes of compact quantum groups.

The modern impetus in the theory of quantum groups came as a result of 
the discovery of new examples of Hopf algebras 
by Drinfeld-Jimbo on the algebraic side \cite{Dr87a} 
and Woronowicz on the analytic side \cite{Wor87a}. These are deformation 
quantizations of the classical Lie algebras and Lie groups, 
and much of the literature on quantum groups had been 
devoted to this approach to quantum groups.

Starting in his thesis \cite{W93}, the author took a different
direction than the traditional deformation quantization method
by viewing quantum groups as intrinsic
objects and found in a series of papers (including \cite{W96b} in
collaboration with Van Daele) several classes of universal compact quantum
groups that {\em can not be obtained} as deformations of Lie groups
or Lie algebras, the
most important of these are the universal compact quantum groups of
Kac type $A_u(n)$ and their orthogonal counterpart $A_o(n)$
\cite{W95a}, the more general universal compact quantum groups $A_u(Q)$
and their self-conjugate counterpart $B_u(Q)$ \cite{W96b,W97b},
where $Q \in GL(n, {\mathbb C})$, and the quantum automorphism groups
$A_{aut}(B, tr)$  \cite{W98a} of finite dimensional $C^*$-algebras $B$ endowed
with a tracial functional $tr$, including the quantum permutation
groups $A_{aut}(X_n)$ on the space $X_n$ of $n$ points. These objects have been
an international focus of study in the subject of compact quantum groups and 
interest in them continues unabated (cf.
\cite{B96a}-\cite{Bi-R-V06},
\cite{Gos09},
\cite{K-Sp09},
\cite{Sol09}-\cite{Voi11}).

The quantum groups $A_u(Q)$ have the remarkable universal property that
can be used to give following alternative and {\em concrete definition} of 
compact matrix quantum groups that was originally defined by Woronowicz \cite{Wor87b} 
more abstractly: a compact matrix quantum group is a quotient
$A_u(Q)/I$, where $I$ is a Woronowicz $C^*$-ideal in $A_u(Q)$. This
means in geometric language
that every compact matrix quantum group (including compact Lie group) is a quantum
subgroup of $A_u(Q)$ for an appropriate choice of $Q$.
In contrast, Drinfeld-Jimbo quantum groups and other
deformations of Lie groups do not enjoy this property,
but are quantum subgroups of $A_u(Q)$.
Similarly, the quantum groups $B_u(Q)$ have the universal property that
every compact matrix quantum group with self-conjugate fundamental representation
is of the form $B_u(Q)/I$, where $I$ is a Woronowicz $C^*$-ideal in $B_u(Q)$.

Without universal compact quantum groups,
the Drinfeld-Jimbo quantum groups and other quantum groups
obtained by deformation would be the end of the story.
However, the outpouring of papers on the universal quantum
groups in the last few years
(e.g.
\cite{B96a}-\cite{Bi-R-V06},
\cite{Gos09},
\cite{K-Sp09},
\cite{Sol09}-\cite{Voi11})
demonstrates depth of the subject:
despite much work achieved so far, we have only seen
the tip of the iceberg and the story is far from the end.

Although compact quantum groups have a satisfactory axiomatic framework and
Drinfeld-Jimbo quantum groups are their special examples, after
the discovery of these new classes of universal compact quantum groups,
it is a natural program to classify simple compact quantum
groups. This program was initiated in \cite{W09},
where it was shown that all compact quantum groups mentioned above are
simple in generic cases.

The {\em main goals} of the program on simple compact quantum groups are:
(1) construct and classify simple compact quantum groups and
their irreducible representations,
(2) understand the structure of simple 
compact quantum groups and structure of 
compact quantum groups in terms of the simple ones, and
(3) develop new applications of simple compact quantum groups
in other areas of mathematics and physics, such as
quantum symmetries in noncommutative geometry and
algebraic quantum field theory. For goal (1), one would
like to develop a theory of simple compact quantum groups that
parallels the Killing-Cartan theory and the Cartan-Weyl theory for
simple compact Lie groups. For this purpose, one must
first construct all simple compact quantum groups. Though the
work so far provides several infinite classes of examples of these, it
should be pointed out that the construction of simple compact
quantum groups is only at the beginning stage for this task at the
moment, as all the simple compact
quantum groups known so far are almost classical in the sense that
their representation rings are isomorphic to those of ordinary
compact groups and in particular are commutative. The first
examples of simple compact quantum groups that are not almost
classical should be directly related to the universal quantum groups $A_u(Q)$ 
(see the footnote after \probref{simpleA_u(Q)}), 
where $Q \in GL(n, {\mathbb C})$ are positive, $n \geq 2$,
though these quantum groups are not simple themselves.
 The representation ring of $A_u(Q)$ is
highly noncommutative, being roughly the free product of two
copies of the ring of integers, according to Banica \cite{B97}. To
construct other simple compact quantum groups, the most natural idea is to
study quantum automorphism groups 
of appropriate quantum spaces, such as those in the author's papers
\cite{W98a,W99b}, the papers of Banica, Bichon, Goswami and their collaborators
\cite{B05a}-\cite{B-Bi-Co07a}, \cite{Bi03a,Bi04a},
\cite{Gos09}, \cite{Bh09}-\cite{Bh-Gos-Sk11}.
In retrospect, both simple Lie groups
and finite simple groups are automorphism groups.
This suggests viability of this approach to the program.
Natural mathematical and physical structures that
have compact quantum automorphic group symmetries
are compact commutative and noncommutative Riemannian manifolds
in the sense of Connes \cite{Cn,Cn96a}.
Such symmetries should be investigated first.

The following heuristic may indicate the depth of the problem on
classification problem of simple compact quantum groups.
The finite dimensional factors are classified by the
{\em discrete set} of natural numbers while
the classification of infinite dimensional von Neumann factors involves
{\em continuous parameters}. Similarly, simple compact Lie groups are
classified by a {\em discrete set} of Cartan matrices,
but the classification of simple
compact quantum groups involves {\em continuous parameters}.
However, since the algebraic structures of compact quantum groups
are richer and more rigid than those of von Neumann algebras,
the classification of simple compact quantum groups might be more accessible
than the classification of infinite dimensional factors.
Even if a classification of simple compact quantum groups
up to isomorphism is unattainable, just as von Neumann factors are far from being
classified up isomorphism, experience has demonstrated
that the study of universal quantum groups and quantum automorphism groups
is fruitful
(cf. \cite{W95a,W96b,W97b,W99b,W02b,W09},
\cite{B96a}-\cite{Bi-R-V06},
\cite{Gos09},
\cite{K-Sp09},
\cite{Sol09}-\cite{Voi11} and references therein),
and other types of classification theories may also be considered, such as
the classification of easy quantum groups \cite{B-Sp09,B-Curr-Sp11,B-Curr-Sp12,Weber12pr}
and the classification of restricted classes of quantum automorphism groups \cite{B-Bi07b}-\cite{B-Bi-Co07b} 
by Banica {\sl et al.}.

An outline of the paper is as follows.
In \secref{Simple}, we review the notion of compact quantum groups 
and simple compact quantum groups.
In  \secref{SimpleEx}, we review examples of simple compact quantum groups constructed so far in 
\cite{W09}. In \secref{Au(Q)Free}, we discuss problem of the fine structure of $A_u(Q)$
and simple quantum quotient groups from $A_u(Q)$, as well as
quantum subgroups from free products of compact quantum groups.
In  \secref{AC-PF}, we give a list of problems related to 
almost classical compact quantum groups and compact quantum groups with property $F$.
In \secref{Graph} and \secref{Isometry}, we discuss the problem of constructing
simple compact quantum groups from quantum automorphisms of finite graphs and other
quantum subgroups of the quantum permutation groups, easy quantum groups
and quantum isometry groups. In \secref{Graph}, several new quantum groups constructed by 
Banica, Curran and Speicher \cite{B-Curr-Sp10} since \cite{W09} are shown to be simple using 
results in Raum \cite{Raum12} and Weber \cite{Weber12pr} along with results in \cite{W09}.

\bigskip

\section{Compact quantum groups and simple compact quantum groups}
\label{Simple}

We first recall the definition of compact (matrix) quantum groups and then the
notion of simple compact quantum groups.

There are several equivalent definitions of compact (matrix) quantum groups, each
has its own advantages over the others. We briefly describe below another
equivalent definition, which has the advantage of letting the reader
to ``visualize'' all compact quantum groups more concretely. This equivalent definition
is essentially in the literature, but not in the explicit form we describe below.

\medskip
\noindent
{\bf Notation:} For elements $u_{ij}$ ($i,j =1, ... n$)  of a
$C^*$-algebra $A$, we define the following elements in the
$n \times n$ matrix algebra $M_n(A)$ over $A$: 
$u:=(u_{ij})_{i,j =1}^n$,  ${\bar u} :=(u^*_{ij})_{i,j =1}^n$,
 $u^t:=(u_{ji})_{i,j =1}^n$ and   $u^*:= {\bar u}^t$, i.e. $u^* = (u^*_{ji})_{i,j =1}^n$.

\bgdf {\rm (cf.  \cite{W93,W95a,W96b,W97b,W02b})}
{\rm 
The {\bf universal compact matrix quantum groups} 
are defined to be the family of pairs $(A_u(Q), \Delta_u)$,  where
$Q \in GL(n, {\mathbb C})$, $Q > 0$, and $u_{ij}$ ($i,j =1, ... n$) are generators of
the universal $C^*$-algebra $A_u(Q)$ that satisfies the following sets of relations:
\bgeqq
    u^* u = I_n = u u^*, \; \; \;
    u^t Q {\bar u} Q^{-1} = I_n = Q {\bar u } Q^{-1} u^t ,
\ndeqq
and $\Delta_u: A_u(Q) \rightarrow A_u(Q) \otimes A_u(Q)$ is the uniquely defined morphism
such that
\bgeqq
\Delta_u(u_{ij}) = \sum_{k=1}^n u_{ik} \otimes u_{kj}.
\ndeqq
}
\nddf

Note that instead of restricting $Q$ to positive matrices,  one can still define
$(A_u(Q), \Delta_u)$ for any invertible  $Q$. For  such $Q$, one has
the following free product decomposition \cite{W02b},
\bgeqq
A_u(Q) \cong A_u(P_1) \ast A_u(P_2) \ast \cdots \ast A_u(P_k)
\ndeqq
for appropriate positive matrices $P_1, P_2, ..., P_k$ with compatible coproducts
as in \cite{W95a}.

\bgdf {\rm (cf. \cite{Wor87b})}
{\rm A {\bf compact matrix quantum group}
is a triple $(A, \Delta, \pi)$,
where    
$\pi: A_u(Q) \rightarrow A$ and $\Delta: A \rightarrow A \otimes A$
are $C^*$-morphisms such that

{\rm (1)} $\pi$ is surjective, 
and

{\rm (2)} $\Delta \pi =  (\pi \otimes  \pi) \Delta_u$.
}
\nddf

Note that using \cite{W95a}, conditions (1) and (2) above is equivalent to (1)
plus the following  condition:

{\rm (2)}$^\prime$ $\Delta_u(\ker(\pi)) \subset ker(\pi \otimes \pi)$.

\medskip

Hence a compact  matrix quantum group can be defined even more simply as a pair
$(A, \pi)$ satisfying (1) and {\rm (2)}$^\prime$.

\medskip
The $C^*$-algebra $A$ in the definition of a
compact matrix quantum group $(A, \Delta, \pi)$
is called a {\bf finitely generated Woronowicz $C^*$-algebra}.
The morphism $\Delta$ is called the {\bf coproduct} of $A$.

It can be shown that there is a Hopf $*$-algebra structure
$({\mathcal A}_u(Q), \Delta, \varepsilon, S)$ on the dense
$*$-subalgebra ${\mathcal A}_u(Q)$ of $A_u(Q)$ such that
\bgeqq
S(u_{ij}) = u^*_{ji}, \;  \;  \;
\varepsilon(u_{ij}) = \delta_{ij}, \;  \;  \;  i, j = 1, 2, \cdots n.
\ndeqq
\noindent
This $*$-Hopf algebra structure induces Hopf $*$-algebra structure
on the dense $*$-subalgebra ${\mathcal A}$ of $A$ in the above definition
of compact matrix quantum group,  and on the dense $*$-subalgebra ${\mathcal A}$
of $A$ in the definition the compact quantum group below.

As in 2.3 of \cite{W95a}, one can define morphism between compact matrix quantum groups
as opposite of morphisms of finitely generated Woronowicz $C^*$-algebras.
Under these morphisms,  compact matrix quantum groups form a category, though
this category is not closed under inverse limits, which leads to the 
following notion of compact quantum groups using 3.1 of \cite{W95a},

\bgdf
{\rm (cf. \cite{BS93,Wor98a})}
{\rm
A {\bf compact quantum group} 
$(A, \Delta)$ is 
an inverse limit of compact matrix quantum groups
$(A_\lambda, \Delta_\lambda, \pi_{\lambda \lambda'})$.
}
\nddf

We will see in the next few paragraphs that the above notion of compact 
quantum groups is equivalent to the elegant and abstract one in \cite{Wor98a}.

As an inductive limit (instead of inverse limit) of
finitely generated Woronowicz $C^*$-algebras $A_\lambda$,
the $C^*$-algebra $A$ in the definition above is called a {\bf Woronowicz $C^*$-algebra}.
Kernels of morphisms between
Woronowicz $C^*$-algebras are called 
{\bf Woronowicz $C^*$-ideals},
which can also be intrinsically defined as in 2.3 of \cite{W95a}.

\medskip
\noindent
{\bf Remark:}
Intuitively, we think of $A=C(G)$, where $G$ is a compact quantum group, even there might be
no points in ``$G$'' other than the identity. We also use  the
notations $A_G=C(G)$ and $C(G_A) = A$.
As in the literature, by abuse of terminology, $A$ is also called a compact (matrix)
quantum group besides being called a (finitely generated) Woronowicz $C^*$-algebra.

\medskip
As usual, one defines the Haar state/measure. One then establishes the
existence and uniqueness of the Haar state/measure on any compact matrix quantum
group just as in \cite{Wor87b}. Using this and 3.3 of \cite{W95a}, one establishes
the existence and uniqueness of the Haar state/measure on any compact quantum
group $G$. In the following, the Haar state on $A = C(G)$ is denoted by $h_G$ or simply
$h$ if no confusion arises.

Using the Haar state/measure $h$,
the Peter-Weyl theory for all compact quantum groups can be developed
as in \cite{Wor98a} or \cite{MaesVD}. As a result, we see that the above 
definition of compact quantum groups is equivalent to the one in \cite{Wor98a}. 

We use ${\mathcal A}_G$ to denote the dense $*$-subalgebra of $A_G$ consisting of 
coefficients of finite dimensional representations of $G$. As a consequence of 
the Peter-Weyl theory for compact quantum groups,  ${\mathcal A}_G$ is a Hopf $*$-algebra.

We need to recall the notion of normal quantum subgroups \cite{W95a,Wnormal}
to define simple compact quantum groups.

Let $(N, \pi)$ be a {\bf quantum subgroup} of a compact quantum group
$G$, which, as defined  in \cite{W93,W95a}, means that
$\pi: C(G) \longrightarrow C(N)$ is a surjection of $C^*$-algebras
such that
$ (\pi \otimes \pi) \Delta_G =\Delta_N \pi $, where $\Delta_G, \Delta_N$
are coproducts of $C(G)$ and $C(N)$, respectively.

Define
$$C(G/N) := \{ a \in C(G) | (id \otimes \pi) \Delta (a) = a \otimes 1_N \},
\; \; \; $$
$$
C({N \backslash G}) := \{ a \in C(G) | (\pi \otimes id) \Delta (a)
= 1_N \otimes a \}, \; \; \; $$
$\Delta$ being coproduct on $C(G)$, $1_N$ the unit of $C(N)$.

\bgdf
\label{prop-normal}
{\rm (cf. \cite{W95a,W09})}
{\rm
We say $N$ is  {\bf normal} in $G$ if it satisfies one of
the equivalent conditions in the proposition below.
}
\nddf

\bgprop
\label{normal-subgroup}
{\rm (cf. \cite{W09})}
Let $N$ be a quantum subgroup of a compact quantum group $G$.
The following conditions are equivalent:

{\rm (1)}  $ C({N \backslash G})$ is a
Woronowicz $C^*$-subalgebra of $C(G)$.

{\rm (2)}  $ C({G/N})$ is a
Woronowicz $C^*$-subalgebra of $C(G)$.

{\rm (3)}  $ C({G/N}) = C({N \backslash G})$.

{\rm (4)}  For every irreducible representation $u^\lambda$ of $G$, either
$h_N \pi (u^\lambda) = I_{d_\lambda}$ or $h_N \pi (u^\lambda) = 0$,
where $h_N$ is the Haar measure on $N$, $d_\lambda = \text{dim}(u^\lambda)$
and $I_{d_\lambda}$ the $d_\lambda \times d_\lambda$
identity matrix.
\ndprop

Among the above four equivalent formulations of the notion of normal quantum
subgroups, condition (4) is the most convenient for our purposes.
If $G$ is a compact group, then the above definition coincides with the
usual notion of closed normal subgroups.

To avoid complications with classification of finite quantum groups,
we want to restrict the notion of simple quantum groups to 
quantum groups that are connected.

\bgdf
\label{simple}
{\rm (cf. \cite{W09})}
{\rm
A compact quantum group $G$ is called {\bf connected}
if for each non-trivial irreducible representation
$u^\alpha  \in \hat{G}$, the $C^*$-algebra $C^*(u^\alpha_{ij})$
generated by the coefficients of $u^\alpha$ is of infinite dimension.

A compact matrix quantum group $G$ is called {\bf simple}
(resp. {\bf absolutely simple}) if it is
connected and has no non-trivial connected normal quantum subgroups
(resp. non-trivial normal quantum subgroups) and
no non-trivial representations of dimension one.
}
\nddf

Note that compact Lie groups are (commutative) examples of compact
quantum groups. It is easy to show that a compact Lie group
is simple (resp. connected) in the usual sense if and only
if it is  simple (resp. connected) in the sense above.

\bgprob 
\label{simpledefprob}
In the definition of simple compact quantum groups, if one replaces the 
condition ``it has no non-trivial representations of dimension one'' 
with the apparently less stringent condition  ``$C(G) \neq C^*(\Gamma)$ where 
$\Gamma$ is a discrete group'', do we get an equivalent definition?  
\ndprob

Note according to \cite{Wor87a}, the condition $C(G) \neq C^*(\Gamma)$ above 
means that $G$ is a nonabelian compact quantum group. 
For a compact Lie group $G$, this condition 
simply means $G$ is a nonabelian Lie group, and  
the answer to the above question is affirmative by Weyl's dimension formula. 


Just as the notion of simple compact Lie groups excludes the torus groups,
the definition of simple quantum groups above (including the alternative 
one formulated in \probref{simpledefprob}) excludes the compact quantum
groups coming from group $C^*$-algebras $C^*(\Gamma)$
of discrete groups $\Gamma$ (i.e., abelian quantum groups).
This is important
because the classification of discrete groups is out of reach.

\bigskip

\section{Examples of simple compact quantum groups that are 
almost classical and have property $F$}
\label{SimpleEx}

The simple compact quantum groups that have been
constructed so far share many properties common to compact Lie
groups.  We recall two of these properties \cite{W09}.

Just as for compact groups,  the {\bf representation ring} (also called the
{\bf fusion ring}) $R(G)$ of a compact quantum group $G$ is a partially ordered algebra over the
integers ${\mathbb Z}$ generated by the irreducible  characters
of $G$.
The set of characters of $G$ is a semi-ring that defines the order of
$R(G)$.

\bgdf
\label{property-F-etc}
{\rm (cf. \cite{W09})}
{\rm
A compact quantum group $G$ is said to have {\bf property $F$} if each
Woronowicz $C^*$-subalgebra of $C(G)$ is of the form $C({G/N})$ for some
normal quantum subgroup $N$ of $G$.

A compact quantum group is called {\bf almost classical}
if its representation ring $R(G)$
is order isomorphic to the representation ring of a compact group.
}
\nddf

In plain language, a compact quantum $G$ is said to have property
$F$ if its quantum function algebra $C(G)$ behaves exactly as the function
algebras of compact groups with respect to normal subgroups.
Note that compact quantum groups $C^*(\Gamma)$ from dual $\Gamma$ of
a discrete group $\Gamma$ do not have this property unless the discrete group
$\Gamma$ is abelian, in which case the group $C^*$-algebra
$C^*(\Gamma)$ is a genuine function algebra over the
Pontryagin dual  $\hat{\Gamma}$ of the discrete abelian group $\Gamma$.

The notion dual to property $F$ is given by following definition, which captures
the property that compact quantum group $C^*(\Gamma)$ has with respect to
normal quantum subgroups:

\bgdf
{\rm (cf. \cite{W09})}
{\rm
A compact quantum $G$ is said to have {\bf property $FD$}
if each of quantum subgroup of $G$ is normal.
}
\nddf

Proofs of assertions in \cite{W09} concerning properties $F$ and $FD$, 
along with other related properties of compact quantum groups, 
can be found in \cite{Wnormal}. 

\bigskip

\noindent
{\bf  Quantum groups $B_u(Q)$}  (cf. \cite{W93,W95a,W96b,W97b,W02b})
\medskip

Keeping the notation for definition of $A_u(Q)$ in
\secref{Simple}. Let $Q \in GL(n, {\Bbb C})$ be such that $ Q \bar{Q} = \pm I_n$, $n \geq 2$.
$B_u(Q) $ is defined to be the universal $C^*$-algebra with generators $u_{ij}$
($i,j = 1, 2, \cdots, n$) that satisfy that following sets of relations:
\bgeqq
    u^* u = I_n = u u^*, \; \; \;
    u^t Q u Q^{-1} = I_n = Q u Q^{-1} u^t.
\ndeqq

It can be shown that there is a well-defined morphism
\bgeqq
\Delta:  B_u(Q)  \rightarrow B_u(Q) \otimes B_u(Q)
\ndeqq
such that
\bgeqq
\Delta(u_{ij}) = \sum_{k=1}^n u_{ik} \otimes u_{kj}.
\ndeqq
and that
$(B_u(Q), \Delta) $ is a compact matrix quantum group. The
quantum groups $B_u(Q)$ have the universal property that
every compact matrix quantum group with self-conjugate fundamental representation
is of the form $B_u(Q)/I$, where $I$ is a Woronowicz $C^*$-ideal in $B_u(Q)$.

\medskip
\noindent
{\bf Note:}  $B_u(Q)$
is also denoted by $A_o(Q^*)$ by Banica {\sl et al}.

When $Q = I_n$, $B_u(Q)$ is  just $A_o(n)$,
the universal orthogonal quantum group of  Kac type
 introduced in 4.5 of \cite{W95a}.  Banica {\sl et al.} also
invented the notation $O^+(n)$ to signify
 $A_o(n) = C(O^+(n))$. 
 
 In addition, when  
$Q 
= 
\left[
\begin{array}{cc}
             0 & -I_n    \\
             I_n & 0
\end{array}
\right] 
$, $B_u(Q)$ is a quantum symplectic group. This is one of the reasons 
that we use the notation $B_u(Q)$ instead of the notation $A_o(Q)$, 
as the latter only captures the special case  $Q=I_n$. 
Following the notation $O^+(n)$ of  Banica {\sl et al.}, we denote 
the above universal symplectic quantum group by $Sp^+(n)$. 
The quantum symplectic group has also appeared in recent work of 
Bhowmick, D'Andrea, Das and D\c{a}browski 
on quantum gauge symmetries (see \cite{Bh-DA-Da-Da12}).

\medskip

For any invertible $Q \in GL(n, {\Bbb C})$  that
does not satisfy $ Q \bar{Q} = \pm I_n$,
$B_u(Q)$ can be defined by the same relations as above, but

\begin{center}
$B_u(Q) \cong A_u(P_1) \ast A_u(P_2) \ast \cdots \ast A_u(P_k) \ast$

$\hspace{3cm}
\ast \,  B_u(Q_1) \ast B_u(Q_2) \ast \cdots \ast B_u(Q_l).$
\end{center}
for certain $P_i>0$ 
and $Q_j$ 
such that $Q_j \bar{Q}_j$'s are scalars (cf. \cite{W02b}).

We note that both $A_u(Q)$ and $B_u(Q)$
can be alternatively described as quantum automorphism groups
of appropriate spaces (cf. \cite{W99b,Bh-DA-Da-Da12}).

\bigskip

\noindent
{\bf Quantum automorphism groups  $A_{aut}(B, \tau)$}
(cf. \cite{W98a})
\medskip

Let $B$ be a of a finite dimensional $C^*$-algebra and
$\phi$ a functional on  $B$.
In general the quantum automorphism group
$A_{aut}(B, \phi)$ that preserves the system $(B, \phi)$ exists,
which is the universal object in the category of compact
quantum groups acting on the system $(B, \phi)$.
However, as shown in Banica \cite{B99b}, only when $\phi$
is {\em the canonical trace $\tau$} does the quantum group 
$A_{aut}(B, \tau)$ have a representation theory that is
relatively easy to describe,
where the trace $\tau$ on $B$ is called {\em canonical}
if it coincides with the restriction to $B$ of the unique tracial state on
the algebra $L(B)$ of operators with $B$ acting by the GNS
representation associated with the trace $\tau$.
If one identifies $B$ with $\bigoplus_{k=1}^m M_{n_k} ({\Bbb C})$,
then
\bgeqq
\tau(\sum_{k=1}^m b_k) = \sum_{k=1}^m \frac{n_k^2}{n} Tr (b_k),
\ndeqq
where $b_k \in M_{n_k}$,  $n$ is the dimension of $B$, 
and $Tr$ is ordinary trace on $M_{n_k}$, i.e.
$Tr (b_k)$ is the sum of the diagonal entries of the matrix $b_k$.

In either case, explicit description of  $A_{aut}(B, \phi)$ or
$A_{aut}(B, \tau)$
in terms of generators and relation is complicated for a general
finite dimensional $C^*$-algebra $B$.
However, when $B=C(X_n)$ is the commutative $C^*$-algebra of functions on the space
where $X_n$ is the space of $n$ points,
the {\bf quantum permutation group} $ A_{aut}(X_n) := A_{aut}(C(X_n))$ has a
surprisingly simple description  in terms of generators and relations:
The $C^*$-algebra $A_{aut}(X_n)$ is generated by self-adjoint
projections $a_{ij}$ such that each row and
column of the matrix $(a_{ij})_{i,j=1}^n$ adds up to 1,  i.e.,
\bgeqq
\begin{aligned}
 & a_{ij}^2 = a_{ij} = a_{ij}^*, \; \; \; i,j = 1, \cdots, n, \\
 & \sum_{j = 1}^{n} a_{ij} = 1 , \; \; \; i = 1, \cdots, n, \\
 & \sum_{i = 1}^{n} a_{ij} = 1 , \; \; \; j = 1, \cdots, n.
\end{aligned}
\ndeqq

Banica {\sl et al.} invented the very convenient geometric notation $S_n^+$
so that $A_{aut}(X_n) = C(S_n^+)$.
In \cite{Bi08a}, Bichon generalizes the quantum permutation groups to
purely algebraic context and establishes the universal property
of the quantum permutation groups in complete generality.

\bgthm
\label{simpleAut}
{\rm (cf. 4.1 and 4.7 in \cite{W09})}

{\rm (a)}
For $Q \in GL(n, {\mathbb C})$ such that $Q \bar{Q} = \pm I_n$ and $n \geq 2$,
$B_u(Q)$ is
almost classical simple compact quantum group with property $F$.

{\rm (b)} Let $B$ be a finite dimensional $C^*$-algebra
endowed with its canonical trace $\tau$ and $\dim(B) \geq 4$, 
then quantum group $A_{aut}(B, \tau)$ is an
almost classical absolutely simple compact quantum groups with property $F$.

\ndthm

As special cases of $B_u(Q)$ and $A_{aut}(B, \tau)$, we have 

\bgcor
{\rm (a)} The universal orthogonal quantum groups 
$O^+(n)$ (for $n \geq 2$) and the universal symplectic quantum groups $Sp^+(n)$ 
(for $n \geq 1$) are simple. 

{\rm (b)} The quantum permutation groups $S_n^+$ are 
simple for $n \geq 4$. 
\ndcor

\medskip

The main ideas used in the proof of \thmref{simpleAut} include

(1) Banica's fundamental work of on the structure of  fusion rings
of these quantum groups (cf. Th\'{e}or\`{e}me 1 and Theorem 4.1 in \cite{B96a,B99b} respectively);

(2) Correspondence between Hopf $*$-ideals and
Woronowicz $C^*$-ideals (cf. 4.2 - 4.3 in \cite{W09});

(3) Reconstruction of a normal quantum group from the
identity in the quotient quantum group (cf. 4.4 in \cite{W09}).

\medskip

In addition to the fusion rings of compact quantum groups such as those considered in  \cite{B96a,B99b}, 
the matters in (2) and (3) above are of interest in their own right
and worth further investigation.
The correspondence between Hopf $*$-ideals and
Woronowicz $C^*$-ideals and related matters in (2) relate
algebraic and analytical aspects of compact quantum groups.
In purely Hopf algebras context,
reconstruction of a normal quantum group from the quotient quantum group 
in (3) has been an issue since 1970's and is related to several other important
and old open questions \cite{Wnormal}.

\bigskip
\noindent
{\bf Quantum groups $K_q$,  $K_q^u$ and $K_J$} (cf. \cite{SV,Soib1,LS,R93b,W96a} and \cite{Rosso87a,Rosso90a})
\medskip

A unified study of the compact quantum groups $K_q$,  $K_q^u$ is due to
Soibelman \& Vaksman, Levendorskii \cite{SV,Soib1,LS}.
The $*$-Hopf algebras ${\mathcal A}_{K_q}$ are algebras of
``representative functions'' of Drinfeld-Jimbo
quantum groups $U_q({\mathfrak g})$ and define in a sense their
``compact real form''.  The quantum group $K_q$
is a deformation of the Poisson Lie group $K(1, 0)$ (cf. \cite{LS}). 
The $*$-Hopf algebras  ${\mathcal A}_{K_q^u}$ are twisting of the
$*$-Hopf algebras ${\mathcal A}_{K_q}$ by an element
$u 
\in \wedge^2 \mathfrak{h}_{\mathbb{R}}$.
The quantum group $K_q^u$ is a deformation of the Poisson Lie group $K(1, u)$ 
(cf. \cite{LS}).  
As shown in \cite{W96a}, $K_q^u$ is an example of Rieffel's deformation from 
action of finite dimensional vector space as conjectured by Rieffel 
(cf. \cite{R93a,R93b} for the background).

Rieffel's quantum group deformation $K_J$ \cite{R93b} depends on $J = S \oplus (-S)$,  
where $S$ is a skew symmetric operator on the Lie algebra (viewed as ${\mathbb{R} }^n$) 
of a torus subgroup of the compact Lie group $K$. 
For appropriate choice of $S$, $K_J$
is a deformation of Poisson Lie group $K(0, u)$ \cite{LS,R93b}.
An action of ${\mathbb{R} }^d : = {\mathbb{R} }^n \times {\mathbb{R} }^n$
on $A=C(K)$ can be constructed and Rieffel's theory 
of deformation for action of ${\mathbb{R} }^d$  \cite{R93a}   
can be applied to obtain $A_J$ \cite{R93b}, also denoted $C(K_J)$.

Precise description of $K_q$,  $K_q^u$ and $K_J$  require more space
than appropriate in this paper. For our purposes, these quantum groups
can be roughly described as follows:

\medskip

(1) The associated dense Hopf $*$-algebras 
${\mathcal A}_{K_q}$,  ${\mathcal A}_{K_q^u}$ and ${\mathcal A}_{K_J}$ 
are the same vector space as the un-deformed/un-twisted ones 
${\mathcal A}_{K}$,  ${\mathcal A}_{K_q}$ and ${\mathcal A}_{K}$ respectively, 
 but the algebras ${\mathcal A}_{K_q}$,  ${\mathcal A}_{K_q^u}$ and ${\mathcal A}_{K_J}$ 
 have deformed products;

(2) The Hopf $*$-algebras ${\mathcal A}_{K_q}$,  ${\mathcal A}_{K_q^u}$ and ${\mathcal A}_{K_J}$ 
have the same coproduct as the un-deformed/un-twisted ones 
${\mathcal A}_{K}$,  ${\mathcal A}_{K_q}$ and ${\mathcal A}_{K}$ respectively;

(3) Representation theories of the deformed/twisted quantum groups $K_q$,  $K_q^u$ and $K_J$ 
are the same as the un-deformed/un-twisted ones $K$,  $K_q$ and $K$  respectively.

\bgthm
\label{simpleK}
{\rm cf. (5.1, 5.4 and 5.6 in \cite{W09})}
If $K$ is a simple compact Lie group,
then $K_q$,  $K_q^u$, $K_J$ are almost classical simple compact quantum groups with property $F$.
\ndthm

The main ideas used in the proof of \thmref{simpleK} include

(1) Representation theory of these quantum groups;

(2) Correspondence between Hopf $*$-ideals and
Woronowicz $C^*$-ideals, as in the proof of \thmref{simpleAut};

(3) Reconstruction of a normal quantum group from the
identity in the quotient quantum group, as in the proof of \thmref{simpleAut};

(4) The normal subgroups of the undeformed Lie group
remain to be normal subgroups of the deformed quantum groups
and explicit identification of normal quantum subgroups of
the deformed quantum groups.

Other deformations of compact Lie groups, though constructed not as systematic as 
the ones considered by Drinfeld-Jimbo, Soibelman {\sl et al.} and Rieffel, 
are scattered in the literature. We believe the general ideas used in the 
proof of \thmref{simpleK} can also be applied to such deformations. 

\bigskip

\section{Quotient quantum groups from $A_u(Q)$ and free products}
\label{Au(Q)Free}

In \cite{W02b}, the quantum groups
$A_u(Q)$ are classified up to isomorphism for positive matrices $Q>0$ and
the quantum groups $B_u(Q)$ are classified up to isomorphism for matrices $Q$ with
$ Q \bar{Q} = \pm I_{n}$.
It is shown that the corresponding $A_u(Q)$ (resp. $B_u(Q)$) is not a free product,
or a tensor product, or a crossed product.
However, for general non-singular matrices $Q$, we have the
decomposition theorem expressing $A_u(Q)$ and $B_u(Q)$ in terms
of free product of the forgoing quantum groups (cf. Theorem 3.1 and 3.3 in \cite{W02b}).
In the light of these results, the following problem seems to be fundamental:

\bgprob
\label{simpleA_u(Q)}
Study further the fine structure of $A_u(Q)$ for
positive matrices $Q \in GL(n, {\mathbb C})$ and $n \geq 2$; 
Determine their simple quotient quantum groups.
\ndprob

A solution of this problem will also provide
the first examples\footnote{After this paper was accepted for publication, 
Alexandru Chirv\u{a}situ informed the author that in the preprint 
``Free unitary quantum groups are (almost) simple'', he showed that 
the quantum group generated by $u_{ij}u^*_{kl}$ and $u^*_{kl}u_{ij}$ is a 
simple compact quantum group with noncommutative representation ring and 
without property $F$, where 
$u_{ij}$ are the generators of $A_u(Q)$ with $Q>0$. 
In the same preprint, he also showed that the quotient quantum group of 
$A_u(Q)$  by its  central subgroup $\mathbb T$ in the proof of 
4.5 in \cite{W09} has no normal quantum subgroups but is not finitely generated, and all 
normal subgroups are subgroups of $\mathbb T$, giving a complete classification of 
simple quotient groups of $A_u(Q)$. 
See also footnote before \probref{ex-ac-pF} below.} of simple compact quantum groups that
are not almost classical (see \secref{SimpleEx}). Note
that the $A_u(Q)$'s have the $1$ dimensional diagonal
 torus $\mathbb T$  as their (connected) normal quantum subgroup, as observed by Bichon
(private communication, cf. 4.5 in \cite{W09}), so they are not simple.
However, they are very close to being simple. For example, they have
no non-trivial irreducible representations of dimension one \cite{B97,W02b}.

It is worth noting that in \probref{simpleA_u(Q)}, simple quotient quantum
groups of $A_u(Q)$ should be easier to determine than
simple quantum subgroups of $A_u(Q)$, since the latter is tantamount to
finding all simple quantum groups due to the universal property of $A_u(Q)$,
which include all simple compact quantum groups in \cite{W09} as reviewed in
\secref{SimpleEx}, simple compact Lie groups, as well as all 
the other unknown simple compact quantum groups.

By investigating the fine structure of concrete quantum groups
such as $A_u(Q)$, one can expect to gain insights
into the structure of general compact quantum groups and simple quantum groups.
Sections \secref{Graph} and \secref{Isometry} below contain more directions
of research on this approach to quantum groups.

Suitable modifications of the method for the proofs
of the main results in section 4 of the paper \cite{W09} should 
yield a solution to \probref{simpleA_u(Q)}. The extra work needed for
this problem that do not appear in \cite{W09} is that
there are more Woronowicz subalgebras in $A_u(Q)$ to consider than therein.
Some preliminary computations of these subalgebras
give optimism to a positive solution of the problem.
One of the main ingredients in this calculation is Banica's 
fundamental result on fusion rules of the irreducible representations
of the quantum group $A_u(Q)$ \cite{B97}.

The general $A_u(Q)$ (resp. $B_u(Q)$) for arbitrary $Q \in GL(n, {\mathbb C})$ 
is not simple if $C^*({\mathbb Z})$ appears in its 
free product decompositions as described in \cite{W02b}.
This is because of the fact that $A_u(Q')=C^*({\mathbb Z})=C({\mathbb T})$ 
(resp. $B_u(Q')=C^*({\mathbb Z}/2{\mathbb Z})$) for $Q' \in GL(1, {\mathbb C})$ and
the following result (\cite{Wnormal}):

\bgprop
\label{normalfreeproduct}
Let $G_1, G_2$ be compact quantum groups.
Let $G= G_1 \hat{*} G_2$ be the  free product compact quantum group
{\rm \cite{W95a}}
underlying
$A_{G_1} * A_{G_2}$.
Let $\pi_1$
be the natural embedding of $G_1$ into $G$
defined by the surjection
$$\pi_1: A_{G_1} * A_{G_2} \rightarrow A_{G_1}, \; \; \; 
\pi_1 = id_1 * \epsilon_2 .  $$
If $G_1$ has at least one irreducible representation
of dimension greater than one,
then $(G_1, \pi_1)$ is {\bf not} a normal quantum subgroup of
$G_1 \hat{*} G_2$.
Otherwise,  $(G_1, \pi_1)$ is
normal in $G_1 \hat{*} G_2$.
\ndprop

The hat in the symbol $\hat{*}$ above
signifies the ``Fourier transform'' of free product $*$
remiscent of classical case in which $G_k = \hat{\Gamma}_k$, where
$\Gamma_k$ are discrete abelian groups.

Note that a compact Lie group being simple means roughly that
it is not a direct product of proper connected subgroups.
A similar result also holds for quantum groups:
If $G_A$ is a simple compact quantum group,
then $A_G$ is not a tensor product
(i.e. $G$ is not a direct product of its non-trivial quantum
subgroups) \cite{W95b,Wnormal}.

However, the proposition above says that the evident quantum subgroups
$(G_1, \pi_1)$ and $(G_1, \pi_2)$ of $G_1 \hat{*} G_2$ are not normal
in $G_1 \hat{*} G_2$ when $G_1$ and
$G_2$ have no non-trivial representations of dimension one.
Along with
other results of the author, this may suggest
that the following problem on the structure of simple quantum groups
has a positive solution.

\bgprob
\label{free-product}
Let $G_1$ and $G_2$ be simple compact quantum groups.
Is $G_1 \hat{*} G_2$ also simple?
\ndprob

The results in author's paper \cite{W95a} should be
useful for a solution of this problem. In particular, Theorem 1.1 there
should play a role, as it did in \cite{B96a,B97,Bi04a,W02b}. 
Note that the formula for the Haar measure on
$G_1 \hat{*} G_2$ and the classification its irreducible representations
are given in explicit formulas in Theorem 1.1 in \cite{W95a}.
According to the postulates in the definition of a normal  quantum subgroup 
(\dfref{prop-normal}), these are important
ingredients in determining whether $G_1 \hat{*} G_2$
has normal quantum subgroups. The results and methods
of the paper \cite{W09} (especially section 4 therein)
should also be useful for this problem.

A positive solution to  \probref{free-product} would have
the following implication: $G_1 \hat{*} G_2$ would be a
{\em simple compact quantum group} when both $G_1$ and $G_2$
are merely {\em simple compact Lie groups}.

The following easier variation of \probref{free-product} is also
of interest and is related to the problem of determining
whether several families of easy quantum groups are simple
(cf. \secref{Graph} below).

\bgprob
\label{free-product2}
Let $G_1$ be a simple compact quantum group and
$G_2$ 
a finite quantum group.
Is $G_1 \hat{*} G_2$ also simple?
\ndprob

Note that in view of the problems in \secref{Graph} below,
it would be interesting to solve \probref{free-product2} for $G_2$ a finite group. 
Also related problems can be formulated for Bichon's free wreath product
of compact quantum groups  \cite{Bi04a} (cf. \secref{Graph} below): 

\bgprob
\label{free-product3}
Let $G_1$ be a simple compact quantum group and
$G_2 = S_n^+$ 
the quantum permutation group.
Is the free wreath product $G_1 \hat{*} S_n^+$ also simple?
\ndprob

\bigskip

\section{Problems related to almost classical compact quantum groups 
and property $F$
}
\label{AC-PF}

As reviewed in \secref{SimpleEx}, the simple compact quantum groups
known so far are almost classical and have property $F$.  Such quantum groups
seems to be most accessible at the moment.
The following problem evidently is  less difficult than the
general problem of classifying simple compact quantum groups and
should be attempted first:

\bgprob
\label{prob-ac-pF}
Classify simple compact quantum groups that are almost classical
and have property $F$.
\ndprob

The following closely related problems should be considered also:

\bgprob
\label{prob-ac-pF2}
{\rm (a)} Classify almost classical simple compact quantum groups.

{\rm (b)} Classify simple compact quantum groups with property $F$.
\ndprob

\bgprob
\label{prob-property-FD}
Does simple compact quantum groups with property $FD$ exist?
If so, construct and classify them.
\ndprob

In \probref{prob-ac-pF} and \probref{prob-ac-pF2},
it would be interesting enough to restrict consideration to
almost classical simple compact quantum groups that
have the same representation rings as simple compact Lie groups.

In another direction, for the apparently more difficult problem
of classifying simple compact quantum groups that
are not almost classical or without property $F$, we do not
have a single example of them. Therefore the following is 
a basic problem\footnote{Alexandru Chirv\u{a}situ informed the author 
that he has since solved \probref{ex-ac-pF} and \probref{prob-noncomm-R} - see 
footnote after \probref{simpleA_u(Q)}.}:

\bgprob
\label{ex-ac-pF}
{\rm (a)}  Construct an example of simple compact quantum group that is not almost classical.

{\rm (b)} Construct an example of simple compact quantum group that does not have property $F$.
\ndprob

A more concrete problem than \probref{ex-ac-pF} is the following

\bgprob
\label{prob-noncomm-R}
Construct simple compact quantum groups with noncommutative representation ring.
\ndprob

For the problems in this section, results on general structure of compact
quantum groups such as those in the previous sections should also be useful.
In this direction, we have the following result (cf. \cite{Wnormal}).

\bgth
Let $G$ be a compact quantum group with property $F$.
Then its quantum subgroups and quotient groups $G/N$
by normal quantum subgroups $N$ also have property $F$.
\ndth

It would be of interest to develop other general results on the structure of compact quantum groups.

\bigskip

\section{Quantum automorphism groups of finite graphs
 and easy quantum groups}
\label{Graph}

How and where do we find quantum groups satisfying the properties in
the problems in \secref{AC-PF} above?  The most natural approach, in our opinion,
is by considering quantum automorphism groups of appropriate commutative
and noncommutative spaces and their quantum subgroups.
Much of the recent work on compact quantum groups falls into this
category. We would like to mention three classes of  these quantum groups:
quantum automorphism groups of finite graphs and other quantum subgroups of
the quantum permutation groups,  easy quantum groups
and quantum isometry groups. The last of these three is discussed in the next section.
We briefly look at the first two in this section.

In part to understand quantum subgroups of
the quantum permutation groups \cite{W98a}, Bichon \cite{Bi03a}
constructed the {\em quantum automorphism groups of finite graphs},
which are quantum subgroups of the former preserving the edges of the graphs.
To further understanding this new class of quantum groups,
Bichon \cite{Bi04a} also defined free wreath product of compact
quantum groups using the quantum permutation groups and
 proved the beautiful formula stating that
the free wreath product of the quantum automorphism group of a graph by
the quantum permutation group $A_{aut}(X_n)$ is the quantum
automorphism group of $n$ disjoint copies of the graph.
Banica also independently studied quantum subgroups of the
quantum permutation groups \cite{B05a,B05b}.
These works lead to a great deal of further studies of quantum automorphism of
finite graphs and quantum subgroups of the quantum permutation groups,
cf. \cite{B05a} - \cite{B-Bi-Co07b} and references therein.

It is instructive to see an immediate application of the quantum
automorphism groups of finite graphs to related problems in the last section.
Clearly, a quantum quotient group $G/N$ of an almost classical quantum
group $G$ is also almost classical. However a {\em crucial} observation is
that a quantum subgroup of an almost classical quantum group needs
not be almost classical. For example,
the quantum permutation groups $A_{aut}(X_n)$ are almost classical
(cf. \cite{B99b,W98a,W09} and \secref{SimpleEx}),
but according to of Bichon \cite{Bi04a}, their quantum subgroups
$A_2({\mathbb Z}/m{\mathbb Z})$, as the quantum authomorphism group 
of certain graphs, 
are not almost classical if $m \geq 3$ (see Corollary 2.7 and
the paragraph following Corollary 4.3 of \cite{Bi04a}).
Though the quantum groups $A_2({\mathbb Z}/m{\mathbb Z})$
are not simple (see proposition 2.6 of \cite{Bi04a}), and
noting that the quantum automorphism group of the trivial graph
(i.e. the quantum permutation group) is simple, it is natural to expect
that it is possible to obtain simple compact quantum groups
that are not almost classical by considering quantum automorphism
groups of other appropriate finite graphs, including other free wreath
products, thus solving \probref{ex-ac-pF}.
Note that since the free wreath is constructed from the free
product, such problems are related to
see \probref{simpleA_u(Q)} and the discussions following it.

In another important and related new direction that has origins in 
works on free wreath product and quantum subgroups of the 
quantum permutation group discussed above, 
Banica and Speicher \cite{B-Sp09} initiated
the study of  easy quantum groups and found several interesting families
of new compact quantum groups. 
A  compact matrix quantum group $G$ with fundamental representation $u$
is called {\bf easy} if

 (1) $G$ lies between 
 $S_n$  and  $O^+(n)$ (cf. \secref{SimpleEx} on $B_u(Q)$ with $Q=I_n$); and

(2) For any $k, l \geq 0$,
$Hom(u^{\otimes k}, u^{\otimes l})$ is linearly generated by operators $T_p$
canonically associated with partitions $p$ of  $k + l$.

See
\cite{B-Sp09} or  \cite{B-Curr-Sp10} or \cite{Weber12pr}
for a description of $T_p$.

In addition to the six families of  free easy quantum groups (also called orthogonal quantum groups)
in \cite{B-Sp09,B-Curr-Sp10},  Weber also found another new family ${B^{\prime}_n}^+$ of
free  easy quantum groups in \cite{Weber12pr}, where we follow his different notation
from \cite{B-Sp09}.   In \cite{Raum12}, Raum  computed the fusion rings of
several easy quantum groups in \cite{B-Sp09} using free products. In works of tour de force,
Banica and Vergnioux computed the fusion rings of
the quantum reflection group $H_n^s$ in  \cite{B-Ver09} and 
the half-librated orthogonal quantum group $O^*(n)$ in  \cite{B-Ver10}.
It would be interesting to see if any of these quantum groups provide solutions
to some of the problems in this section.
As the fusion rings of  $H_n^s$ and  $O^*(n)$ are
non-commutative, the following problems seem to be most appealing:

\bgprob
{\rm (a)} Is the quantum reflection group $H^s_n$  simple?

{\rm (b)} Is the half-librated orthogonal quantum group $O^*(n)$ simple?
\ndprob

A positive answer to either (a) or (b) would provide the first simple
compact quantum group that is not almost classical and
has noncommutative representation ring.

In the light of Theorem 4.1 in Raum \cite{Raum12} and 
3.1 and 3.2 in Weber \cite{Weber12pr} whose notation we follow,
the relevant quantum groups considered in that theorem are either not simple or
rely on solution of problems in \secref{Au(Q)Free}.
For instance, using their results above along with our \thmref{simpleAut}, we see 
that 

(1) $B_n^+$ is simple for $n \geq 3$; 

(2) ${B'_n}^+$ ($n \geq 3$) and ${S'_n}^+$ ($n \geq 4$) 
are simple modulo two components (i.e. disconnected); 

(3) ${B^{\#+}_n}$ is not simple because of \propref{normalfreeproduct}, 
but it seems not far from being simple (a concept to be precised) because of 
4.1.(3) in \cite{Raum12} and 3.2.(a) in \cite{Weber12pr}.

It is not clear if the quantum groups
 $H^+_n$, $H^*_n$,  ${H^{(s)}_n}$ and  ${H^{[s]}_n}$  are simple.
Note that as a free wreath product, the quantum group  $H^+_n$ is related to 
the general question \probref{free-product3}.

\bigskip

\section{Quantum isometry groups of commutative  and noncommutative Riemannian spaces}
\label{Isometry}

After Banica's initial investigation of quantum isometry groups of finite spaces \cite{B05a,B05b}, a
recent conceptual breakthrough in compact quantum groups is Debashish Goswami's  \cite{Gos09} theory of
quantum isometry groups of spectral triples \`a la Connes \cite{Cn}, where the universal quantum
groups $A_u(Q)$ plays an essential r\^{o}le in the proof of existence.
Using this notion he computed with Bhowmick \cite{Bh-Gos09a,Bh09} several examples of 
quantum isometry groups and found that they are either isomorphic to classical isometry groups or
are among the examples studied earlier by Rieffel and the author  \cite{R93b, W96a}.
Using the same circle of ideas they subsequently developed \cite{Bh-Gos09b}
an improved notion of quantum isometry group without relying on existence of a good 
Laplacian as required in \cite{Gos09}.  

The quantum isometry groups computed so far are either classical
groups, or known quantum groups,  or combinations of
both based on free product or tensor product.
Because of this, Goswami made a rigidity conjecture earlier stating 
to the effect that connected spaces do not admit non-trivial quantum symmetries.
On the other hand, Huichi Huang \cite{Huang12pr} has shown that
connected non-smooth metric spaces admit faithful action even by
quantum permutation groups, disproving this earlier conjecture. 
Most striking of all is that Goswami, along with his collaborators, 
have recently shown in a series of papers \cite{Gos11-12pr,DGosJ12pr} 
a modified rigidity result stating that a connected and oriented Riemannian manifold
does not have quantum symmetries other than the classical ones.

As a fundamental new concept, one would naturally
wonder if fundamentally new compact quantum groups
can be constructed using quantum isometry groups.
In the light of results of Huang and Goswami {\sl et al.} above,
one should look at quantum isometry groups of 
non-smooth metric spaces or disconnected Riemannian manifolds 
that are non-classical (i.e. not compact  groups).
It would be interesting to find if any such
quantum groups are simple:

\bgprob
\label{simpleQI}
Construct examples of simple quantum isometry groups of
non-smooth metric spaces or disconnected Riemannian manifolds.
\ndprob

It is conceivable that quantum isometry groups of disconnected 
Riemannian manifolds will be related to quantum subgroups of 
the quantum permutation group as considered in \secref{Graph}, 
since quantum permutation group can permute the connected components 
in quantum manner just as it does on the finite space.  

Another direction of research in the theory of quantum isometry group is the following. 
In the newly developed notion of quantum isometry group in \cite{Bh-Gos09b},
the quantum groups in the categories ${\bf Q}'_R$ and ${\bf Q}'$
that are used to define the quantum isometry group there 
does not carry a  $C^*$-algebraic action on the spectral triple, and as 
a result, the universal object (i.e, the quantum isometry group) does not 
carry a $C^*$-algebraic action in general. See \cite{Bh-Gos10b} for an 
example of this  situation. 
(We refer the reader to the above cited papers for 
detailed description of ${\bf Q}'_R$ and ${\bf Q}'$  due to space limitation.) 

To address this problem, we believe the categories ${\bf Q}'_R$ and ${\bf Q}'$ 
are too large, and propose the following alternative for the notion of 
quantum isometry groups 
that will always carry $C^*$-algebraic action. 

First, by a {\bf compact quantum transformation group} $(A, \alpha, u)$
of a compact type spectral triple $({\mathcal B}, {\mathcal H}, D)$ we mean
 a compact quantum transformation group (cf. \cite{W98a}) $(A, \alpha)$ of $B$
that satisfies

(QT1) There is a unitary representation $u \in L({\mathcal H} \otimes A_G)$
of $G_A$ on $\mathcal H$ such that $(\rho, u)$ is $\alpha$-covariant:
$(\rho \otimes 1) \alpha (b) = u (\rho(b) \otimes 1) u^*, \; \; b \in B ;$

(QT2) $D$ is an intertwiner of $u$ with itself:
$(D \otimes 1) u = u (D \otimes 1).$

A {\bf morphism} from another quantum transformation
group $(\tilde{A}, \tilde{\alpha}, \tilde{u})$
to $(A, \alpha, u)$ is defined to be a morphism
$\pi$ from $(\tilde{A}, \tilde{\alpha})$ to $(A, \alpha)$ (cf. \cite{W98a})
that satisfies
$\tilde{u} = (id_{\mathcal H} \otimes \pi) u.$

The above defines {\bf category $\bf {\mathcal C}$ of compact quantum transformation groups}
of $({\mathcal B}, {\mathcal H}, D)$ with objects $\{ (A, \alpha, u) \}$ and morphisms $\{ \pi \}$.
Similar to \cite{W98a}, the {\bf quantum isometry group} of
$({\mathcal B}, {\mathcal H}, D)$
is defined to be a universal object of category $\bf {\mathcal C}$ if it exists.

As in \cite{W98a}, universal object does not always exist in $\bf {\mathcal C}$ in general.
Therefore, as in \cite{W98a} and \cite{Bh-Gos09b}, consider {\bf measured spectral triple}
$({\mathcal B}, {\mathcal H}, D, \phi)$
with a (usually positive) functional $\phi$ on ${\mathcal B}$ and consider the category
$\bf {\mathcal C}_\phi$ of compact quantum transformation groups that satisfies
(QT1)-(QT2) {\em and}

(QT3) $(\phi \otimes  {\rm id} ) (\alpha(b))= \phi(b)1_A$ for $b \in {\mathcal B}$

The {\bf quantum isometry group} of
$({\mathcal B}, {\mathcal H}, D, \phi)$
is defined to be the universal object of category $\bf {\mathcal C}_\phi$ if it exists.
The following can be proved and justifies in part
the above notion quantum isometry group.

\bgprop
{\rm (1)} Let $({\mathcal B}, {\mathcal H}, D)$
be the spectral triple associated with a {\bf compact} Riemannian manifold
$M$. Then the universal object in the category of compact transformation
{\bf groups} of $({\mathcal B}, {\mathcal H}, D)$
in the sense above 
is the {\bf isometry group} of $M$.

{\rm (2)} For an arbitrary spectral triple $({\mathcal B}, {\mathcal H}, D)$,
the universal object in the category of compact transformation
{\bf groups} of $({\mathcal B}, {\mathcal H}, D)$ in the sense above
is the compact {\bf group} of isometries of $({\mathcal B}, {\mathcal H}, D)$
in the sense of Connes (see p6200 of {\rm \cite{Cn95a}}).
\ndprop

Since composition of two continuous maps is continuous, 
the quantum isometry group defined above is evidently contained in the
quantum isometry groups defined by Goswami \cite{Gos09}
and Bhowmick and Goswami \cite{Bh-Gos09b}.
Because of this,  and using ideas in \cite{W98a} and \cite{Bh-Gos09b}, 
the following seems to be quite plausible:

\bgprob
\label{quantum isometry groups of noncommutative spaces}
For a measured spectral triple $({\mathcal B}, {\mathcal H}, D, \phi)$,
show that universal object exists in the category $\bf {\mathcal C}_\phi$.
\ndprob

Many other problems naturally arise:

(i) Calculate the quantum isometry groups in the sense above for
the classical spaces such as the flat spheres, tori, and other Riemannian spaces;

(ii) Calculate the quantum isometry groups in the sense above for
the noncommutative tori;

(iii) Study the properties of quantum isometry groups in general and apply
them to study other properties of noncommutative spaces;

(iv) Construct simple compact quantum groups through the study of quantum isometry groups.

As in the definitions of quantum isometry groups by Goswami and Bhowmick, 
and because a quantum isometry group in our sense is contained in theirs, 
examples of simple quantum isometry groups in our sense might need to be constructed 
out of non-smooth metric spaces or disconnected Riemannian manifolds, as in 
\probref{simpleQI}.

\bigskip

\subsection*{Acknowledgments}
The author would like to thank Hanfeng Li and Yi-Jun Yao for invitations
to Chongqing University and Fudan University respectively during the summer of 2012,
where the author had the opportunity to work on this paper.
The author would also like to thank Hanfeng Li and Dechao Zheng at
Chongqing University, and Yi-Jun Yao and Guoliang Yu at Fudan University for
their hospitality and for creating a congenial atmosphere.
He also wish to record his thanks to Huichi Huang, Pan Ma and
Qinggang Ren for making his visit enjoyable.
The author is indebted to Professor Marek Bozejko for drawing
his attention to easy quantum groups during the conference in September 2011 celebrating
Professor Woronowicz's seventieth birthday.


\begin{thebibliography}{99}



\bibitem{BS93} Baaj, S. and Skandalis, G.:
{\rm Unitaires multiplicatifs et dualit\'e pour les produits
crois\'es de $C^*$-alg\`ebres,}
{\em Ann. Sci. Ec. Norm. Sup.} {\bf 26} (1993), 425-488.


\bibitem{B96a} Banica, T.:
{\rm Th\'eorie des repr\'esentations du groupe quantique
compact libre $O(n)$,}
{\em C. R. Acad. Sci. Paris} t. {\bf 322}, Serie I (1996), 241-244.


\bibitem{B97} Banica, T.:
{\rm Le groupe quantique compact libre $U(n)$,}
{\em Commun. Math. Phys.} {\bf 190} (1997), 143-172.

\bibitem{B99b} Banica, T.:
Symmetries of a generic coaction.
{\em Math. Ann. } {\bf 314} (1999), 763-780.

\bibitem{B05a} Banica, T.:
Quantum automorphism groups of homogeneous graphs.
{\em J. Funct. Anal.} {\bf 224} (2005), no. 2, 243--280.


\bibitem{B05b}
Banica, T.:
Quantum automorphism groups of small metric spaces.
{\em Pacific J. Math.} {\bf 219} (2005), no. 1, 27--51.


\bibitem{B-Bi07a}
Banica, T.,  Bichon, J.:
Free product formulae for quantum permutation groups.
{\em J. Inst. Math. Jussieu} {\bf 6} (2007), no. 3, 381--414.


\bibitem{B-Bi07b}
Banica, T.,  Bichon, J.:
Quantum automorphism groups of vertex-transitive graphs of order $\leq11$.
{\em J. Algebraic Combin.} {\bf 26} (2007), no. 1, 83--105.


\bibitem{B-Bi07pr}
Banica, T.,  Bichon, J.:
Quantum groups acting on 4 points,
{\em J. Reine Angew. Math.}, {\bf 626} (2009), 74-114.

\bibitem{B-Bi-Ch07}
Banica, T.,  Bichon, J.,  Chenevier, G.:
Graphs having no quantum symmetry.
{\em Ann. Inst. Fourier (Grenoble)} {\bf 57} (2007), no. 3, 955--971.


\bibitem{B-Bi-Co07a}
Banica, T.,  Bichon, J.,  Collins, B.:
The hyperoctahedral quantum group.
{\em J. Ramanujan Math. Soc.} {\bf  22}  (2007),  no. 4, 345--384.


\bibitem{B-Bi-Co07b}
Banica, T.,  Bichon, J.,  Collins, B.:
Quantum permutation groups: a survey,
{\em Banach Center Publications} 78, 13-34, 2007.



\bibitem{B-Bi-Nat12}
Banica, T.,  Bichon, J.,   Natale, S.:
Finite quantum groups and quantum permutation groups.
{\em Adv. Math.} {\bf  229} (2012), 3320-3338


\bibitem{B-Co07a}
Banica, T.,  Collins, B.:
Integration over quantum permutation groups.
{\em J. Funct. Anal.} {\bf 242} (2007), no. 2, 641--657.


\bibitem{B-Co07b}
Banica, T.,  Collins, B.:
Integration over compact quantum groups.
{\em Publ. Res. Inst. Math. Sci.} {\bf 43} (2007), no. 2, 277--302.


\bibitem{B-Co08}
Banica, T.,  Collins, B.:
Integration over the Pauli quantum group.
{\em J. Geom. Phys.} {\bf  58} (2008), 942-961.




\bibitem{B-Curr-Sp10}
Banica, T.,  Curran, S., Speicher, R.:
Classification results for easy quantum groups.
{\em Pacific J. Math.} {\bf  247} (2010), no. 1, 1-26.


\bibitem{B-Curr-Sp11}
Banica, T.,  Curran, S., Speicher, R.:
Stochastic aspects of easy quantum groups.
{\em Probab. Theory Related Fields} {\bf  149} (2011), 435-462.




\bibitem{B-Curr-Sp12}
Banica, T.,  Curran, S., Speicher, R.:
De Finetti theorems for easy quantum groups.
{\em Ann. Probab.} {\bf  40} (2012), 401-435.



\bibitem{B-Mo07}
Banica, T.,  Moroianu, S.:
On the structure of quantum permutation groups.
{\em Proc. Amer. Math. Soc.} {\bf 135} (2007), no. 1, 21--29.


\bibitem{B-Sp09}
Banica, T., Speicher, R.:
Liberation of orthogonal Lie groups.
{\em Adv. Math.} {\bf 222} (2009), 1461-1501.


\bibitem{B-Ver09}
Banica, T. and Vergnioux, R.:
Fusion rules for quantum reflection groups.
{\em J. Noncommut. Geom.} {\bf 3}  (2009), 327-359.

\bibitem{B-Ver10}
Banica, T. and Vergnioux, R.:
 Invariants of the half-liberated orthogonal group.
 {\em Ann. Inst. Fourier} {\bf 60}  (2010), 2137-2164.


\bibitem{Bh09}
Bhowmick, J.:
Quantum Isometry Group of the n tori.
{\em Proc. Amer. Math. Soc.} {\bf 137}  (2009), no. 9, 3155-3161.


\bibitem{Bh-DA-Da11}
Bhowmick, J., D'Andrea, F., D\c{a}browski, L.:
Quantum isometries of the finite noncommutative geometry of the standard model.
{\em Comm. Math. Phys.} {\bf 307} (2011), no. 1, 101-131

\bibitem{Bh-DA-Da-Da12}
Bhowmick, J., D'Andrea, F.,  Das,  B., D\c{a}browski, L.:
Quantum gauge symmetries in Noncommutative Geometry. 
arXiv:1112.3622 


\bibitem{Bh-Gos09a}
Bhowmick, J. and Goswami, D.:
Quantum Isometry Groups: Examples and Computations,
{\it Commun.Math.Phys.} {\bf 285} (2009), 421-444,



\bibitem{Bh-Gos09b}
Bhowmick, J. and Goswami, D.:
Quantum group of orientation-preserving Riemannian isometries.
{\em J. Funct. Anal.} {\bf 257}  (2009), no. 8, 2530-2572.


\bibitem{Bh-Gos10a}
Bhowmick, J. and Goswami, D.:
Quantum isometry groups of the Podles spheres.
{\em J. Funct. Anal.} {\bf 258}  (2010), no. 9, 2937-2960.


\bibitem{Bh-Gos10b}
Bhowmick, J., Goswami, D.:
 Some counterexamples in the theory of quantum isometry groups.
 {\em Lett. Math. Phys.} {\bf 93}  (2010), no. 3, 279-293.



\bibitem{Bh-Gos-Sk11}
Bhowmick, J., Goswami, D., Skalski, A.:
Quantum Isometry Groups of 0- Dimensional Manifolds.
{\em Trans. Amer. Math. Soc.} {\bf 363}  (2011), no. 2, 901-921.

\bibitem{Bi03a} Bichon, J.:
Quantum automorphism groups of finite graphs.
{\em Proc. Amer. Math. Soc.} {\bf 131} (2003), no. 3, 665--673



\bibitem{Bi04a} Bichon, J.:
Free wreath product by the quantum permutation group.
{\em Algebr. Represent. Theory} {\bf 7} (2004),  343--362.

\bibitem{Bi08a} Bichon, J.:
Algebraic quantum permutation groups.
{\em Asian-Eur. J. Math.} {\bf 1}  (2008),  no. 1, 1--13.


\bibitem{Bi-R-V06}
Bichon, J.,  De Rijdt, A., Vaes, S.:
Ergodic coactions with large multiplicity and monoidal equivalence of quantum groups.
{\em Comm. Math. Phys.} {\bf 262} (2006), 703--728.

\bibitem{Cn} Connes, A.:
{\em Noncommutative Geometry,} Academic Press, 1994.


\bibitem{Cn95a} Connes, A.:
{\rm Noncommutative geometry and reality,}
{\em J. Math. Phys.} {\bf 36}:11 (1995), 6194-6231


\bibitem{Cn96a} Connes, A.:
{\rm Gravity coupled with matter and the foundation of
non commutative geometry,}
{\em Commun. Math. Phys.} {\bf 182} (1996), 155-176.

\bibitem{DGosJ12pr}
Das, B. and Goswami, D and  Joardar, S.:
Rigidity of action of compact quantum groups II.
arXiv:1206.1718 	


\bibitem{Dr87a} Drinfeld, V.G.:
{Quantum groups},
{\em in} {Proc. of the ICM-1986, Berkeley,} Vol I,
Amer. Math. Soc., Providence, R.I., 1987, pp798--820.


\bibitem{EnSch75a} Enock, M. et Schwartz, J.-M.:
{\rm Une dualit\'e dans les alg\`ebres de von Neumann,}
{\em Bull. Soc. Math. France. Suppl. M\'emoire} {\bf 44} (1975), 1-144.

\bibitem{EnSch92a} Enock, M. et Schwartz, J.-M.:
{\em Kac algebras and duality of locally compact groups.} (English summary)
With a preface by Alain Connes. With a postface by Adrian Ocneanu. 
Springer-Verlag, Berlin, 1992. x+257 pp. ISBN: 3-540-54745-2


\bibitem{Fred89a} Fredenhagen, K. and Rehren, K.H. and Schroer, B.:
{Superselection sectors with braid group statistics and exchange algebras I,}
{\em Comm.  Math. Phys.} {\bf 125} (1989), 201-226.


\bibitem{Gos09}
Goswami, D.:
Quantum  Group of Isometries in Classical and Noncommutative Geometry,
{\em Comm. Math. Phys.} {\bf 285} (2009), no. 1, 141-160.



\bibitem{Gos11-12pr}
Goswami, D.:
 Rigidity of action of compact quantum groups I, III
 	arXiv:1106.5107, arXiv:1207.6470


	
	
\bibitem{Huang12pr}
Huang, Huichi:
 Faithful compact quantum group actions on connected compact metrizable spaces,
 arXiv:1202.1175

\bibitem{Kac63a65a} Kac, G.:
{\rm Ring groups and the duality principle, I, II}
{\em Proc. Moscow Math. Soc.} {\bf 12} (1963), 259-303;
ibid. {\bf 13} (1965), 84-113.


\bibitem{KV74a} Kac, G. and Vainerman, L. I.:
{\rm Nonunimodular ring groups and Hopf Von Neumann algebras,}
{\em Math. USSR Sb.} {\bf 23} (1974), 185-214.


\bibitem{K-Sp09} K\"ostler, C., Speicher, R.:
A noncommutative de Finetti theorem: Invariance under quantum
permutations is equivalent to freeness with amalgamation.
{\em Comm. Math. Phys.} {\bf 291} (2009), no. 2, 473-490.

\bibitem{KusVaes00a} Kustermans, Johan and Vaes, Stefaan
{\rm Locally compact quantum groups.}
{\em  Ann. Sci. \'{E}ole Norm. Sup.} (4) {\bf 33} (2000), no. 6, 837-934.



\bibitem{Lev1}  Levendorskii, S.:
Twisted algebra of functions on compact quantum group and
their representations,
{\em Algebra i analiz} {\bf 3}:2 (1991), 180-198.
{\em St. Petersburg Math. J.} {\bf 3}:2 (1992), 405-423.


\bibitem{LS}  Levendorskii, S. and Soibelman, Y.:
{\rm Algebra of functions on compact quantum groups,
Schubert cells, and quantum tori,}
{\em Commun. Math. Phys.} {\bf 139} (1991), 141-170.

\bibitem{MaesVD}
Maes, Annand Van Daele, Alfons:
{\rm Notes on compact quantum groups.}
{\em Nieuw Arch. Wisk.} (4) {\bf 16} (1998), no. 1-2, 73-112.

\bibitem{MaNaWor03a}  Masuda, T.; Nakagami, Y.; Woronowicz, S. L.
{\rm A $C^*$-algebraic framework for quantum groups.}
{\em Internat. J. Math.} {\bf 14} (2003), no. 9, 903-1001

\bibitem{PW90} Podles, P. and Woronowicz S. L.:
{\rm Quantum deformation of Lorentz group,}
{\em Commun. Math. Phys.} {\bf 130} (1990), 381-431.


\bibitem{Raum12}
Raum, S.;
Isomorphisms and fusion rules of orthogonal free quantum groups and their free complexifications.
{\em Proc. Amer. Math. Soc.} {\bf 140} (2012), no. 9, 3207-3218.


\bibitem{R93a}  Rieffel, M.:
{\rm Deformation quantization for actions of ${\mbox{\bf R}}^n$,}
{\em Memoirs A.M.S.} no. {\bf 506}, 1993.




\bibitem{R93b}  Rieffel, M.:
{\rm Compact quantum groups associated with toral subgroups,}
{\em Contemp. Math.} {\bf 145} (1993), 465-491.


\bibitem{Rosso87a}  Rosso, M.:
{\rm Comparaison des groupes $SU(2)$ quantiques de Drinfeld et Woronowicz,}
{\em C. R. Acad. Sci. Paris} {\bf 304} (1987), 323-326.


\bibitem{Rosso90a}  Rosso, M.:
{\rm Alg\`{e}bres enveloppantes quantifi\'{e}es, groupes quantiques
compacts de matrices et calcul differentiel non-commutatif,}
{\em Duke Math. J.} {\bf 61} (1990), 11-40.

\bibitem{Soib1}  Soibelman, Y.:
{Algebra of functions on compact quantum group and its representations},
{\em Algebra i analiz} {\bf 2}:1 (1990), 190-212.
{\em Leningrad Math. J.} {\bf 2}:1 (1991), 161-178.


\bibitem{SV}   Soibelman, Y. and Vaksman, L.:
{\rm The algebra of functions on quantum $SU(2)$,}
{\em Funct. Anal. ego Pril.} {\bf  22}:3 (1988), 1-14.
{\em Funct. Anal. Appl.} {\bf 22}:3 (1988), 170-181.


\bibitem{Sol09} Soltan, P.M.:
 Quantum families of maps and quantum semigroups on finite quantum spaces.
 J. Geom. Phys. 59 (2009), no. 3, 354-368.

\bibitem{Va-Ven08}
 Vaes, S. and Vennet, N.V.:
Identification of the Poisson and Martin boundaries of orthogonal discrete quantum groups.
{\em Journal of the Institute of Mathematics of Jussieu} {\bf 7} (2008), 391-412.

\bibitem{Va-Ven10}
 Vaes, S. and Vennet, N.V.:
Poisson boundary of the discrete quantum group $\widehat{A_u(F)}$.
{\em Compos. Math.} {\bf 146} (2010), no. 4, 1073--1095.


\bibitem{Va-Ver07}
Vaes, S. and Vergnioux, R.:
The boundary of universal discrete quantum groups, exactness and factoriality.
{\em Duke Mathematical Journal} {\bf 140} (2007), 35--84.

\bibitem{Va-Ver11pr}
Vergnioux, R.,  Voigt , C.:
The K-theory of free quantum groups.
arXiv:1112.3291

\bibitem{Voi11}
 Voigt, C.:
The Baum-Connes conjecture for free orthogonal quantum groups.
{\em Adv. Math.} {\bf 227}  (2011), 1873 - 1913.



\bibitem{W96b} Van Daele, A. and Wang, S.Z.:
Universal quantum groups,
{\it International J. of Math.} {\bf 7} no2 (1996), 255-264.


\bibitem{W93} Wang, S. Z.:
{\em General Constructions of Compact Quantum Groups,}
Ph.D. Thesis, University of California at Berkeley, March, 1993.


\bibitem{W95a}   Wang, S.Z.:
Free products of compact quantum groups,
{\it Commun. Math. Phys.} {\bf 167} (1995), 671-692.


\bibitem{W95b} Wang, S.Z.:
Tensor products and crossed products of compact
quantum groups,  {\it Proc. London Math. Soc.} {\bf 71} (1995), 695-720.


\bibitem{W96a} Wang, S.Z.:
Deformations of compact quantum groups via Rieffel
quantization, {\it Commun. Math. Phys.} {\bf 178}:3 (1996), 747-764.

\bibitem{W97b} Wang, S.Z.:
{\rm Problems in the theory of quantum groups,}
in {\it Quantum Groups and Quantum Spaces},
Banach Center Publication 40 (1997), Inst. of Math., Polish Acad. Sci.,
Editors: R. Budzynski, W. Pusz, and S. Zakrzewski. pp67-78



\bibitem{W98a} Wang, S.Z.:
Quantum symmetry groups of finite spaces,
{\it Commun. Math. Phys.} {\bf 195}:1 (1998), 195-211.

\bibitem{W99b} Wang, S.Z.:
Ergodic actions of universal quantum groups
on operator algebras,
{\it Commun. Math. Phys.} {\bf 203} (1999), 481-498.


\bibitem{W99c} Wang, S.Z.:
Classification of quantum groups $SU_q(n)$,
{\it J. London Math. Soc.} {\bf 59}:2 (1999), 669-680

\bibitem{W02b}  Wang, S.Z.:
Structure and isomorphism classification of quantum
groups $A_u(Q)$ and $B_u(Q)$,
{\it Journal of Operator Theory} {\bf 48}(3) (2002), 573-583



\bibitem{W09} Wang, S.Z.:
Simple compact quantum groups I,
{\em J. Funct. Anal.} {\bf 256}  (2009), no. 10, 3313-3341.

\bibitem{Wnormal} Wang, S.Z.:
Equivalent notions of normal quantum subgroups,
compact quantum groups with properties $F$ and $FD$,  and other applications,
preprint. 


\bibitem{Weber12pr}
Weber, M.:
 On the classification of easy quantum groups - The nonhyperoctahedral and the half-liberated case.
 arXiv:1201.4723


\bibitem{Wor87a}  Woronowicz, S. L.:
{\rm Twisted $SU(2)$ group. An example of noncommutative
differential calculus,}
{\em Publ. RIMS, Kyoto Univ.} {\bf 23} (1987), 117-181.


\bibitem{Wor87b}  Woronowicz, S. L.:
{\rm Compact matrix pseudogroups,}
{\em Commun. Math. Phys.} {\bf 111} (1987), 613-665.



\bibitem{Wor88a}  Woronowicz, S. L.:
{\rm Tannaka-Krein duality for compact matrix pseudogroups.
Twisted $SU(N)$ groups,}
{\em Invent. Math.} {\bf 93} (1988), 35-76.



\bibitem{Wor98a} Woronowicz, S. L.:
{\rm Compact quantum groups,}
in {\em Quantum Symmetries},
Les Houches Summer School-1995, Session LXIV,
A. Connes, K. Gawedzki, J. Zinn-Justin, ed.,
Elsevier Science,
Amsterdam, 1998; pp845-884.



\end{thebibliography}
\end{document}